%%%%%%%%%%%%%%%%%%%%%%%%%%%%%%%%%%%%%%%%%%%%%%%%%%%%%%%%%%%%%%%%
%        T. Holm,  R. Kessar, M. Linckelmann
%        
%        14-11-a4tilde.tex
%        
%        last modified: 14.11.2005
%%%%%%%%%%%%%%%%%%%%%%%%%%%%%%%%%%%%%%%%%%%%%%%%%%%%%%%%%%%%%%%%

\input amstex
%     amsdefinitions.tex
%     __________________

%                 Akzente

%  Kalligraphische Buchstaben  -  Bold Letters

\define\CO{\Cal O}

\define\Zed{\Bbb Z}

                         % Tensorproducts

\define\tenO{\underset {\Cal O}\to\otimes}

                         % Groups and group algebras

\define\OG{{\Cal O G}}
\define\OGb{{\Cal O Gb}}

           % Charaktere

                         % Subscripts on the right

                         % Subscripts on the left

                         % Subscripts on both sides

                         % Categories

                         %Functors and operators

\define\Id{{\operatorname{Id}}}

\define\Proj{{\operatorname{Proj}}}

                         %Rings and groups

                    % Operatornames

\define\Irr{{\operatorname{Irr}\nolimits}}

% dim  is already defined

% ker is already defined

\input xypic

\define\OAfourtilde{{\CO \tilde A_4}}

\countdef\bibliono=3\bibliono=0
\def\entry{\advance\bibliono by 1}
\entry\edef\Brouea{{\number\bibliono}}
\entry\edef\Broueb{{\number\bibliono}}
\entry\edef\BrPu{{\number\bibliono}}
\entry\edef\CaPi{{\number\bibliono}}
\entry\edef\CaPicorr{{\number\bibliono}}
\entry\edef\Erd{{\number\bibliono}}
\entry\edef\Feit{{\number\bibliono}}
\entry\edef\Puig{{\number\bibliono}}
\entry\edef\Thev{{\number\bibliono}}

\documentstyle{amsppt}
\magnification=1200
\tolerance=10000
\hbadness=10000  
\pageheight{192mm}
\pagewidth{132mm}

\topmatter  
\title  Blocks with a quaternion defect group over a $2$-adic ring:
the case $\tilde A_4$ 
\endtitle
\rightheadtext{Blocks with a quaternion defect group}
\author Thorsten Holm, Radha Kessar, Markus Linckelmann  \endauthor

\address             \vbox{\hbox{Thorsten Holm}
                      \hbox{\hskip 3mm Department of Pure Mathematics}
                      \hbox{\hskip 3mm University of Leeds}
                      \hbox{\hskip 3mm Leeds, LS2 9JT}
                      \hbox{U.K.}}
\endaddress

\address             \vbox{\hbox{Radha Kessar}
                      \hbox{\hskip 3mm Department of Mathematical Sciences}
                      \hbox{\hskip 3mm Meston Building}
                      \hbox{\hskip 3mm Aberdeen, AB24 3UE}
                      \hbox{U.K.}}
\endaddress

\address             \vbox{\hbox{Markus Linckelmann}
                      \hbox{\hskip 3mm Department of Mathematical Sciences}
                      \hbox{\hskip 3mm Meston Building}
                      \hbox{\hskip 3mm Aberdeen, AB24 3UE}
                      \hbox{U.K.}}
\endaddress

\abstract
Except for blocks with a cyclic or Klein four defect group, it is not known
in general
whether the Morita equivalence class of a block algebra over a field of
prime characteristic determines that of the corresponding block algebra
over a $p$-adic ring. We prove this to be the case when the defect group
is quaternion of order $8$ and the block algebra over an algebraically closed
field $k$ of characteristic $2$ is Morita equivalent to $k\tilde A_4$.
The main ingredients are Erdmann's classification of tame blocks [\Erd]
and work of Cabanes and Picaronny [\CaPi, \CaPicorr] on perfect isometries
between tame blocks.
\endabstract

\endtopmatter  

\document

\heading Introduction \endheading

\bigskip

Throughout these notes, $\CO$ is a complete discrete valuation ring
with algebraically closed residue field $k$ of characteristic $2$
and with quotient field $K$ of characteristic $0$.
According to Erdmann's classification in [\Erd], if $G$ is a finite group
and if $b$ is a block of $\OG$ having the quaternion group $Q_8$
of order $8$ as defect group, then the block algebra $kG\bar b$
is Morita equivalent to either $kQ_8$ or $k\tilde A_4$ or
the principal block algebra of $k\tilde A_5$, where here $\bar b$ is
the canonical image of $b$ in $kG$. In the first case the block is
{it nilpotent} (cf. [\BrPu]), and it
follows from Puig's structure theorem of nilpotent blocks in [\Puig] that
$\OG b$ is Morita equivalent to $\CO Q_8$. In the remaining two 
cases one should expect that $\OG b$ is Morita equivalent to
$\CO \tilde A_4$ or the principal block algebra of $\CO\tilde A_5$,
respectively. We show this to be true in one of these two cases under
the assumption that $K$ is large enough:

\bigskip

\proclaim{Theorem A} Let $G$ be a finite group, and let $b$ be a block
of $\OG$ having a quaternion defect group of order $8$. Denote by
$\bar b$ the image of $b$ in $kG$. Assume that $KGb$ is split.
If $kG\bar b$ is Morita equivalent to $k\tilde A_4$ then $\OGb$ is
Morita equivalent to $\OAfourtilde$.
\endproclaim

\bigskip

By Cabanes-Picaronny [\CaPi, \CaPicorr], 
in the situation of Theorem A there is a perfect
isometry between the character groups of $\OGb$ and of $\OAfourtilde$. 
Thus Theorem A is a consequence of the following slightly
more general Theorem which characterises $\OG b$ in terms of 
its center, its character group and $k\tilde A_4$; see the end of
this section for more details regarding the notation.

\bigskip

\proclaim{Theorem B} Let $A$ be an $\CO$-free $\CO$-algebra such that
$K\tenO A$ is split semi-simple and such that $k\tenO A$ is Morita
equivalent to $k\tilde A_4$. Assume that there is 
an isometry $\Phi : \Zed\Irr_K(A)\cong \Zed\Irr_K(\OAfourtilde)$ which
maps $\Proj(A)$ to $\Proj(\OAfourtilde)$ such that the map sending
$e(\chi)$ to $e(\Phi(\chi))$ for every $\chi \in \Irr_K(A)$ induces
an $\CO$-algebra isomorphism of the centers $Z(A)\cong Z(\OAfourtilde)$.
Then $A$ is Morita equivalent to $\OAfourtilde$.
\endproclaim

\bigskip

Theorem B is in turn a consequence of the more precise Theorem C,
describing $A$ in terms of generators and relations:

\bigskip

\proclaim{Theorem C} 
Let $A$ be a basic $\CO$-free $\CO$-algebra such that
$K\tenO A$ is split semi-simple and such that $k\tenO A$ is 
isomorphic to $k\tilde A_4$. Assume that there is 
an isometry $\Phi : \Zed\Irr_K(A)\cong \Zed\Irr_K(\OAfourtilde)$ which
maps $\Proj(A)$ to $\Proj(\OAfourtilde)$ such that the map sending
$e(\chi)$ to $e(\Phi(\chi))$ for every $\chi \in \Irr_K(A)$ induces
an $\CO$-algebra isomorphism of the centers $Z(A)\cong Z(\OAfourtilde)$.
Then     $A$ is  isomorphic to the unitary  $\CO$-algebra   with set of
generators  
$\{e_0, e_1, e_2, \beta, \gamma, \delta, \eta, 
\lambda, \kappa\}$    of $A$, such that $e_0, e_1, e_2$ are 
pairwise orthogonal
idempotents whose sum is $1$ and satisfying  the following relations:
$$\beta = e_0\beta = \beta e_1,\ \gamma = e_1\gamma = \gamma e_0\ ;$$
$$\delta = e_1\delta = \delta e_2,\ \eta = e_2\eta = \eta e_1\ ;$$
$$\lambda = e_2\lambda = \lambda e_0, \ \kappa = e_0\kappa = \kappa e_2\ ;$$
$$\beta\delta = -2\kappa + \kappa\lambda\kappa\ ;\ \ 
\eta\gamma = -2\lambda + \lambda\kappa\lambda\ ;\ \ 
\delta\lambda = -2\gamma + \gamma\beta\gamma\ ;$$
$$\kappa\eta = -2\beta + \beta\gamma\beta\ ;\ \  
\lambda\beta = -2\eta + \eta\delta\eta\ ;\ \ 
\gamma\kappa = -2\delta + \delta\eta\delta\ ;$$
$$\gamma\beta\delta = -4\delta + 2\delta\eta\delta\ ;\ \ 
\delta \eta \gamma = -4\gamma + 2\gamma \beta \gamma \ ;\  \
\lambda \kappa \eta = -4\eta + 2\eta\delta \eta \ ;$$
$$\beta\gamma \kappa = -4\kappa + 2\kappa\lambda\kappa\ ;\ \ 
\eta\delta\lambda = -4\lambda + 2\lambda \kappa \lambda\ ;\  \ 
\kappa\lambda\beta = -4\beta + 2\beta\gamma\beta\ ; $$
$$\eta\gamma\beta = -4\eta + 2\eta\delta\eta\ ;\ \ 
\beta \delta \eta  = -4\beta + 2 \beta \gamma \beta \ ;\  \
\delta \lambda \kappa  = -4\delta + 2\delta \eta \delta\ ;$$
$$\lambda \beta\gamma  = -4\lambda + 2\lambda\kappa \lambda\ ;\ \ 
\kappa\eta\delta = -4\kappa + 2 \kappa \lambda \kappa\ ;\  \ 
\gamma \kappa\lambda = -4\gamma + 2\gamma\beta \gamma\ ; $$
$$\beta\delta\lambda \beta = -8\beta + 4\beta\gamma\beta \ ;\ \ 
\delta\lambda \beta\delta = -8\delta + 4\delta \eta \delta \ ;\  \ 
\lambda\beta\delta\lambda = -8\lambda + 4\lambda\kappa \lambda \ ; $$

\endproclaim

\bigskip

When reduced modulo $2$, these relations  seem  to be more than those 
occuring
in Erdmann's work [\Erd] over $k$ (we recall these more precisely in 
$\S 2$);  
but  they are not, since the 
extra relations   over $k$ can be deduced  from those given by 
Erdmann.  We need to add in extra relations over $\CO$ in order to 
ensure that  the algebra we construct is  $\CO$-free of the right rank.

Since $\CO\tilde A_4$ fulfills the hypotheses of Theorem C it
follows that $A \cong \CO\tilde A_4$, hence Theorem C indeed implies 
Theorem B.
The proof of Theorem C is given  at the end of Section 2.

\bigskip

\noindent
{\bf Notation.} If $A$ is an $\CO$-algebra such that $K\tenO A$ is split
semi-simple, we denote by $\Irr_K(A)$ the set of characters of the
simple $K\tenO A$-modules, viewed as central functions from $A$ to $\CO$
and we denote by $\Irr_k(k\tenO A)$ the set of isomorphism classes of
simple $k\tenO A$-modules.
We denote by $\Zed\Irr_K(A)$ the group of characters of $A$, and we 
denote by
$\Proj(A)$ the subgroup of $\Zed\Irr_K(A)$ generated by the characters
of the projective indecomposable $A$-modules. We denote by $L^0(A)$ 
the subgroup
of $\Zed\Irr_K(A)$ of all elements which are orthogonal to $\Proj(A)$ 
with
respect to the usual scalar product in $\Zed\Irr_K(A)$.
For any $\chi\in\Irr_K(A)$, we denote by $e(\chi)$ the corresponding
primitive idempotent in $Z(K\tenO A)$. If $A = \OG$ for some finite group
$G$ we have the well-known formula 
$$e(\chi) = \frac{\chi(1)}{\vert G\vert}
\underset{x\in G}\to\sum \ \chi(x^{-1})x\ .$$
We refer to [\Brouea, \Broueb] for the concept and basic properties of
perfect isometries, and to [\Thev] for general block theoretic background
material.

\bigskip

\heading 1 Characters and perfect isometries of $\OAfourtilde$ \endheading

\bigskip

We identify $\tilde A_4 = Q_8\rtimes C_3$. Let $t$ be a generator of $C_3$
and let $y$ be an element of order $4$ in $Q_8$. Set $z = y^2$; that is,
$z$ is the unique central involution of $\tilde A_4$. Then the seven elements
$1$, $z$, $y$, $t$, $t^2$, $tz$, $t^2z$ are a complete set of representatives
of the conjugacy classes in $\tilde A_4$. 

Let $\omega$ be a primitive third root of unity in $\CO$. The character table
of $\tilde A_4$ is as follows:

\bigskip

$$\matrix & & 1 & z & y & t & t^2 & tz & t^2z \\
& & & & & & & & \\
\eta_0 & & 1 & 1 & 1 & 1 & 1 & 1 & 1 \\
\eta_1 & & 1 & 1 & 1 & \omega & \omega^2 & \omega & \omega^2 \\
\eta_2 & & 1 & 1 & 1 & \omega^2 & \omega & \omega^2 & \omega \\
\eta_3 & & 3 & 3 & -1 & 0 & 0 & 0 & 0 \\
\eta_4 & & 2 & -2 & 0 & -\omega^2 & -\omega & \omega^2 & \omega \\
\eta_5 & & 2 & -2 & 0 & -\omega & -\omega^2 & \omega & \omega^2 \\
\eta_6 & & 2 & -2 & 0 & -1 & -1 & 1 & 1 \\
\endmatrix$$

\bigskip

The algebra $\OAfourtilde$ has three simple modules $T_0$, $T_1$, $T_2$,
up to isomorphism. Choosing for $T_0$ the trivial module and after possibly
exchanging the notation for $T_1$, $T_2$, the ordinary decomposition
matrix of $\OAfourtilde$ is as follows:
$$\pmatrix 1 & 0 & 0 \\ 0 & 1 & 0 \\ 0 & 0 & 1 \\ 1 & 1 & 1 \\
 1 & 1 & 0 \\ 1 & 0 & 1 \\ 0 & 1 & 1 \\
 \endpmatrix $$

\bigskip

The Cartan matrix of $\OAfourtilde$ is the product of the decomposition
matrix with its transpose, hence equal to
$$\pmatrix 4 & 2 & 2\\ 2 & 4 & 2\\ 2 & 2 & 4\endpmatrix $$

\bigskip

Let $e_0$, $e_1$, $e_2$ be primitive idempotents in $\OAfourtilde$
such that $\OAfourtilde e_i$ is a projective cover of $T_i$, 
$0 \leq i \leq 2$. By the above decomposition matrix,
the characters of the projective indecomposable $\OAfourtilde$-modules
$\OAfourtilde e_i$ are
$$\eta_0 + \eta_3 + \eta_4 + \eta_5\ ,$$
$$\eta_1 + \eta_3 + \eta_4 + \eta_6\ ,$$
$$\eta_2 + \eta_3 + \eta_5 + \eta_6\ ,$$
respectively. Their norm is $4$, and the differences of any two different
characters of projective indecomposable $\OAfourtilde$-modules yields the 
following
further elements in $\Proj(\OAfourtilde)$ having also norm $4$:
$$\eta_0-\eta_1+\eta_5-\eta_6\ ,$$
$$\eta_0-\eta_2+\eta_4-\eta_6\ ,$$
$$\eta_1-\eta_2+\eta_4-\eta_5\ .$$
It is easy to check, that up to signs, these are all elements in
$\Proj(\OAfourtilde)$ having norm $4$.

\bigskip

\bigskip

A self-isometry $\Phi$ of $\Zed\Irr_K(\OAfourtilde)$ maps every
$\eta_i$ to $\epsilon_i\eta_{\pi(i)}$ for some signs $\epsilon_i
\in \{1,-1\}$ and a permutation $\pi$ of $\{0,1,\ldots,6\}$. In other
words, $\Phi$ is determined by the permutation $\tau$ of the set
$\{1,-1\}\times\{0,1,\ldots,6\}$ satisfying $\tau(1,i) = (\epsilon_i,
\pi(i))$ and $\tau(-1,i) = (-\epsilon_i,\pi(i))$ for all
$i$, $0\leq i\leq 6$. If we write $i$, $-i$ instead of $(1,i)$,
$(-1,i)$, respectively, this becomes $\tau(i) = \epsilon_i\pi(i)$
and $\tau(-i) = -\epsilon_i\pi(i)$, with the usual cancellation
rules for signs. In this way, every self-isometry $\Phi$ of
$\Zed\Irr_K(\OAfourtilde)$ gets identified to a permutation of the
set of symbols $\{i, -i | 0\leq i\leq 6\}$.

\bigskip

A perfect self-isometry of $\Zed\Irr_K(\OAfourtilde)$ 
preserves necessarily $\Proj(\OAfourtilde)$. The next Proposition
implies that the converse is true, too:

\bigskip

\proclaim{Proposition 1.1} The group of all perfect self-isometries
of $\Zed\Irr_K(\OAfourtilde)$ is equal to the group of all self-isometries
of $\Zed\Irr_K(\OAfourtilde)$ 
which preserve $\Proj(\OAfourtilde)$. This group is
generated by $-\Id$ together with the set of permutations
$$(0,1,2)(4,6,5)\ ,$$
$$(1,2)(4,5)\ ,$$
$$(2,-3)(5,-6)\ .$$
\endproclaim

\bigskip

Every algebra automorphism of $\OAfourtilde$
induces a permutation on $\Irr_K(\OAfourtilde)$ which is in fact a
perfect isometry on $\Zed\Irr_K(\OAfourtilde)$. 
Since $\eta_1$ has degree $1$, it is an algebra
homomorphism from $\OAfourtilde$ to $\CO$, and hence the map sending
$x\in \OAfourtilde$ to $\eta_1(x) x$ is an algebra automorphism of
$\OAfourtilde$ whose inverse sends $x\in \OAfourtilde$ to $\eta_2(x) x$.
The following statement is an immediate consequence from the
character table of $\OAfourtilde$:

\bigskip

\proclaim{Lemma 1.2} Let $\gamma$ be the algebra automorphism
of $\OAfourtilde$ defined by $\gamma(x) = \eta_1(x) x$ for all $x\in 
\OAfourtilde$. The permutation $\pi$ of $\{0,1,\ldots,6\}$ defined by
$\eta_i\circ\gamma = \eta_{\pi(i)}$ is equal to $\pi = (0,1,2)(4,6,5)$.
\endproclaim

\bigskip

The anti-automorphism of $\OAfourtilde$ sending $x\in \tilde A_4$ to
$x^{-1}$ induces also a permutation of the set $\Irr_K(\OAfourtilde)$,
and this is also a perfect isometry (this holds for any finite group).
This permutation can also be read off the character table:

\bigskip

\proclaim{Lemma 1.3} Let $\iota$ be the algebra anti-automorphism of 
$\OAfourtilde$ mapping
$x\in \tilde A_4$ to $x^{-1}$. The permutation $\pi$ of 
$\{0,1,\ldots,6\}$
defined by $\eta_i\circ\iota = \eta_{\pi(i)}$ is equal to 
$\pi = (1,2)(4,5)$.
\endproclaim

\bigskip

\demo{Proof of 1.1} 
The first two permutations are perfect
isometries by 2.2 and 2.3, respectively.
An easy but painfully long verification shows that the bicharacter
sending $(g,h)\in\tilde A_4\times\tilde A_4$ to
$$\eta_0(g)\eta_0(h)+\eta_1(g)\eta_1(h)-\eta_2(g)
\eta_3(h)-\eta_3(g)\eta_2(h)+\eta_4(g)\eta_4(h)-
\eta_5(g)\eta_6(h)-\eta_6(g)\eta_5(h)$$
is perfect; that is, its value at any $(g,h)$
is divisible in $\CO$ by the orders of
$C_{\tilde A_4}(g)$ and $C_{\tilde A_4}(h)$ and it vanishes if
exactly one of $g$, $h$ has odd order.
Thus the isometry given by the permutation $(2,-3)(5,-6)$
is perfect. It remains to show that these permutations, together with
$-\Id$,
generate the group of all self-isometries which preserve 
$\Proj(\OAfourtilde)$.

We described above a complete list of all
elements in $\Proj(\OAfourtilde)$ having norm $4$. Since the characters
of the projective indecomposable modules are in that list, a self-isometry
$\Phi$ of $\Zed\Irr_K(\OAfourtilde)$ preserves $\Proj(\OAfourtilde)$ if and
only if it permutes this set of norm $4$ elements. 

Let $\Phi$ be a self-isometry of $\Zed\Irr_K(\OAfourtilde)$ which
preserves $\Proj(\OAfourtilde)$. Then $\Phi$ preserves
also the group $L^0(\OAfourtilde)$ of generalised characters which
are orthogonal to all characters in $\Proj(\OAfourtilde)$. 
Up to signs, the complete
list of elements in $L^0(\OAfourtilde)$ having norm $3$ is
$$\eta_0+\eta_1-\eta_4\ ,\ \eta_0+\eta_2-\eta_5\ ,\ \eta_0-\eta_3+\eta_6\ ,$$
$$\eta_1+\eta_2-\eta_6\ ,\ \eta_1-\eta_3+\eta_5\ ,\ \eta_2-\eta_3+\eta_4\ .$$
Up to signs again, the complete list of elements in $L^0(\OAfourtilde)$
having norm $4$ is
$$\eta_0+\eta_1+\eta_2-\eta_3\ ,$$
$$\eta_0-\eta_1-\eta_5+\eta_6\ ,\ 
\eta_0-\eta_2-\eta_4+\eta_6\ ,\ 
\eta_0+\eta_3-\eta_4-\eta_5\ ,$$ 
$$\eta_1-\eta_2-\eta_4+\eta_5\ , \ 
\eta_1+\eta_3-\eta_4-\eta_6\ ,\ \eta_2+\eta_3-\eta_5-\eta_6\ .$$
The first norm $4$ element in this list, $\eta_0+\eta_1+\eta_2-\eta_3$,
is the only norm $4$ element which is orthogonal to all other norm
$4$ elements in $L^0(\OAfourtilde)$. Thus $\Phi$ has to permute
the characters $\eta_0$, $\eta_1$, $\eta_2$, $\eta_3$ amongst
each other. 

Suppose first that $\Phi$ fixes $\eta_3$. Then, by composing $\Phi$
with a suitable product of powers of the first two permutations
in the statement, we may assume that $\Phi$ fixes $\eta_0$, $\eta_1$,
$\eta_2$ up to signs. By considering the first of the above norm $4$
elements in $L^0(\OAfourtilde)$ we get that $\Phi$ fixes $\eta_0$, $\eta_1$,
$\eta_2$ all with positive signs. By considering the norm $3$ elements in
$L^0(\OAfourtilde)$, it follows that $\Phi$ fixes also $\eta_4$, $\eta_5$
and $\eta_6$ with positive signs. Thus a self-isometry of 
$\Zed\Irr_K(\OAfourtilde)$ which preserves $\Proj(\OAfourtilde)$ and which fixes
$\eta_3$ is in the group generated by the set of two permutations
$(0,1,2)(4,6,5)$ and $(1,2)(4,5)$. 

Suppose next that $\Phi$ does not fix $\eta_3$. By precomposing $\Phi$
with a suitable of power of $(0,1,2)(4,6,5)$ 
we may assume that $\Phi$ sends $\eta_2$ to $-\eta_3$.
By composing $\Phi$ with a suitable power of $(0,1,2)(4,5,6)$
we may assume that $\Phi$ fixes $\eta_0$, up to a sign. 
Since $\Phi$ preserves the norm $4$ element
$\eta_0+\eta_1+\eta_2-\eta_3$, we necessarily have $\Phi(\eta_0)
= \eta_0$. Then $\Phi$ maps $\eta_1$ either to $\eta_1$ or $\eta_2$
(with positive signs, again because of that same norm $4$ element).
In the first case, $\Phi$ fixes both $\eta_0$, $\eta_1$, and by checking
the norm $3$ elements in $L^0(\OAfourtilde)$ one gets $\Phi = (2,-3)(5,-6)$.
In the second case, again checking on norm $3$ elements, one gets
$\Phi = (1,2,-3)(4,5,-6)$, but this is already the product of $(1,2)(4,5)$
and $(2,-3)(5,-6)$.
$\square$ \enddemo

\bigskip

\heading 2 The algebra $A$ \endheading

\bigskip

Let $A$ be a basic $\CO$-algebra fulfilling the hypotheses
of Theorem B; that is, $K\tenO A$ is split semi-simple, 
$k\tenO A$ is isomorphic to $k\tilde A_4$,
and there is an isometry $\Zed\Irr_K(A)\cong\Zed\Irr_K(\OAfourtilde)$
mapping $\Proj(A)$ to $\Proj(\OAfourtilde)$ and inducing an isomorphism
$Z(A)\cong Z(\OAfourtilde)$. There is a ``compatible choice" for these isomorphisms:

\bigskip

\proclaim{Proposition 2.1} There is an algebra isomorphism $\alpha :
k\tenO A\cong k\tilde A_4$ and an isometry $\Phi : \Zed\Irr_K(A)\cong \Zed\Irr_K
(\OAfourtilde)$ mapping $\Proj(A)$ to $\Proj(\OAfourtilde)$
with the following properties:
\smallskip
(i) $\Phi$ maps $\Irr_K(A)$ onto $\Irr_K(\OAfourtilde)$; that is, all signs are
$+1$.
\smallskip
(ii) The map sending $e(\chi)$ to $e(\Phi(\chi))$ for every $\chi\in
\Irr_K(A)$ induces an isomorphism $Z(A)\cong Z(\OAfourtilde)$.
\smallskip
(iii) For any primitive idempotents $e \in A$ and 
$f\in \OAfourtilde$ and every $\chi\in\Irr_K(A)$ such that $\alpha(\bar e)
= \bar f$
we have $\chi(e) = \Phi(\chi)(f)$; that is, $A$ and $\OAfourtilde$ have the
same decomposition matrices through $\alpha$ and $\Phi$.
\endproclaim

\demo{Proof} The $\CO$-rank of $A$ is $24$ and also the sum of the squares
of the seven irreducible $K$-linear characters of $A$; thus every irreducible
character of $A$ has degree smaller than $5$. Also, there is no character
of degree $4$ because $24 - 4^2 = 8$ cannot be written as a sum of six
squares of the six remaining characters. But there must be a character
of degree $3$; if not, $24$ would be the sum of seven squares all either
$1$ or $4$, which is not possible. Thus the squares of the six remaining
characters add up to $24-3^2 = 15$, and the only way to do this is with
three characters of degree $1$ and three characters of degree $2$.

This proves
that the character degrees of the irreducible characters of $A$ and of
$\OAfourtilde$ coincide for some bijection $\Irr_K(A) \cong \Irr_K(\OAfourtilde)$.
Since the decomposition matrix of $A$ multiplied with its transpose
yields the Cartan matrix of $A$ - which is equal to that of $k\tilde A_4$ -
the algebra $A$ has in fact the same decomposition matrix as $\OAfourtilde$
for a suitable bijection $\Phi : \Irr_K(A)\cong \Irr_K(\OAfourtilde)$
and the bijection $\Irr_k(k\tenO A)\cong \Irr_k(k\tilde A_4)$ induced
by $\alpha$. Extend $\Phi$ to a $\Zed$-linear isomorphism
$\Zed\Irr_K(A)\cong\Zed\Irr_K(\OAfourtilde)$, still denoted by $\Phi$.
By construction, $\Phi$ sends the characters of the projective indecomposable
$A$-modules to the characters of the projective indecomposable
$\OAfourtilde$-modules; in particular, $\Phi$ maps $\Proj(A)$ to
$\Proj(\OAfourtilde)$. It remains to see that the map sending $e(\chi)$
to $e(\Phi(\chi))$ for every $\chi\in\Irr_K(A)$ induces an isomorphism
$Z(A)\cong Z(\OAfourtilde)$. For any $i$, $0\leq i\leq 6$, denote by
$\chi_i$ the irreducible character of $A$ such that $\Phi(\chi_i) = \eta_i$.
As in the proof of 1.1, we have a distinguished norm $4$ element
in $L^0(A)$ which is orthogonal to all other norm $4$ elements in
$L^0(A)$, namely $\chi_0+\chi_1+\chi_2-\chi_3$. Thus, if $\Psi :
\Zed\Irr_K(A)\cong \Zed\Irr_K(\OAfourtilde)$ is some isometry mapping
$\Proj(A)$ to $\Proj(\OAfourtilde)$ and inducing an isomorphism
$Z(A)\cong Z(\OAfourtilde)$, then $\Psi(\chi_0+\chi_1+\chi_2-\chi_3)
= \pm(\eta_0+\eta_1+\eta_2-\eta_3)$. By Proposition 1.1, there is
a perfect self-isometry $\mu$ of $\Zed\Irr_K(\OAfourtilde)$
such that $\Phi = \mu\circ\Psi$.
$\square$ \enddemo

\bigskip

\noindent{\bf Remark 2.2.} If we assume that $A$ is Morita equivalent to
some block algebra with $Q_8$ as defect group, then Proposition 2.1
follows also from the work of Cabanes and Picaronny in [\CaPi, \CaPicorr].

\bigskip

Since $k\tenO A \cong k\tilde A_4$, 
the quiver of $A$ is the same as that of $k\tilde A_4$,
thus of the following form:

$$
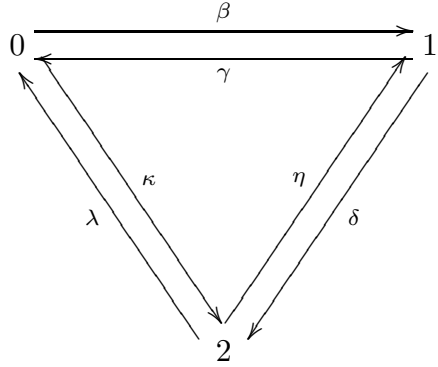
\xymatrix{ 0 \ar@<1ex>[rrrr]^\beta \ar@<1ex>[rrddd]^\kappa& & & & 
1 \ar@<1ex>[llll]^\gamma \ar@<1ex>[llddd]^\delta \\
 & & & & \\
 & & & & \\
 & & 2 \ar@<1ex>[uuull]^\lambda \ar@<1ex>[uuurr]^\eta & & }$$

Write $\bar a$ for the image of $a\in A$ in $\bar A =
k\tenO A \cong k\tilde A_4$. The generators $\beta$, $\gamma$, $\delta$,
$\kappa$, $\lambda$, $\eta$ can be chosen such that their images in
$\bar A$ fulfill the following relations:
$$\bar\beta\bar\delta = \bar\kappa\bar\lambda\bar\kappa\ ,$$
$$\bar\eta\bar\gamma = \bar\lambda\bar\kappa\bar\lambda\ ,$$
$$\bar\delta\bar\lambda = \bar\gamma\bar\beta\bar\gamma\ ,$$
$$\bar\kappa\bar\eta = \bar\beta\bar\gamma\bar\beta\ ,$$
$$\bar\lambda\bar\beta = \bar\eta\bar\delta\bar\eta\ ,$$
$$\bar\gamma\bar\kappa = \bar\delta\bar\eta\bar\delta$$

and

$$\bar\gamma\bar\beta\bar\delta = \bar\delta\bar\eta\bar\gamma =
\bar\lambda\bar\kappa\bar\eta = 0 \ .$$

In order to determine the algebra structure of $A$, we have to ``lift"
these relations over $\CO$.

\bigskip

We fix an algebra isomorphism $\alpha : k\tenO A \cong k\tilde A_4$ and
an isometry $\Phi : \Zed\Irr_K(A)\cong \Zed\Irr_K(\OAfourtilde)$ satisfying
the conclusions of Proposition 2.1. We denote by $\chi_i$ the unique
irreducible $K$-linear character of $A$ such that $\Phi(\chi_i) = \eta_i$
for all $i$, $0\leq i\leq 6$.

The characters $\eta_0$, $\eta_1$, $\eta_2$, $\eta_3$ of $\OAfourtilde$
have height zero, the characters $\eta_4$, $\eta_5$, $\eta_6$ have height
one. Thus, via the isomorphism of the centers induced by $\Phi$,
it follows that for $0\leq i\leq 3$ we have
$8e(\chi_i) \in A$, and for $4\leq j\leq 6$ we have $4e(\chi_j)\in A$.
We can in fact describe an $\CO$-basis of $Z(A)$ in terms of the 
centrally primitive idempotents $e(\chi_i)$.
The strategy is now to play off the descriptions of $Z(k\tenO A)$ in
terms of the generators in the quiver and of $Z(A)$
in terms of the centrally primitive idempotents $e(\chi_i)$.

\bigskip

\proclaim{Lemma 2.3} The following elements of $Z(K\tenO A)$ are all
contained in the radical $J(Z(A))$:
$$s = 2e(\chi_4) + 2e(\chi_5) + 2 e(\chi_6)\ ,$$
$$z_0 = 4e(\chi_2)+4e(\chi_3)+2e(\chi_4)\ ,$$
$$z_1 = 4e(\chi_1)+4e(\chi_3)+2e(\chi_5)\ ,$$
$$z_2 = 4e(\chi_0)+4e(\chi_3)+2e(\chi_6)\ ,$$
$$y_0 = 4e(\chi_1)+4e(\chi_2)+2e(\chi_4)+2e(\chi_5)\ ,$$
$$y_1 = 4e(\chi_0)+4e(\chi_2)+2e(\chi_4)+2e(\chi_6)\ ,$$
$$y_2 = 4e(\chi_0)+4e(\chi_1)+2e(\chi_5)+2e(\chi_6)\ .$$
Moreover, for any two different $i$, $j$ in $\{0,1,2\}$ the set
$$\{1, z_i, z_j, s, 8e(\chi_3), 4e(\chi_{i+4}), 4e(\chi_{j+4})\}$$
is an $\CO$-basis of $Z(A)$.
\endproclaim

\demo{Proof} In view of Proposition 2.1
we may assume that $A = \OAfourtilde$. 
This is just an explicit verification, using the character
table of $\tilde A_4$. One verifies first that $z_0 \in A$.
By symmetry, this implies that $z_1$, $z_2$ are also in $A$.
Then $y_0 = z_0 + z_1 - 8e(\chi_3)$ is in $A$, similarly
for the $y_1$, $y_2$. An equally easy computation shows that
$s\in A$. Thus all the given elements belong to $Z(A)$.
None of these elements is invertible, so they all belong
to $J(Z(A))$ because $Z(A)$ is local.

In order to see the last statement on the basis of $Z(A)$, we 
may assume that $i=0$ and $j=1$.
For any $x\in \tilde A_4$ denote by $\underline x$ the conjugacy
class sum of $x$ in $\CO\tilde A_4$. The orthogonality relations imply the
well-known formula
$$\underline x = \underset{0\leq m\leq 6}\to\sum\ \frac{\chi_m(
\underline x^{-1})}{\chi_m(1)} e(\chi_m)\ .$$
Thus, for the seven conjugacy classes in $\tilde A_4$, we have
$$\underline 1 = e(\chi_0)+e(\chi_1)+e(\chi_2)+e(\chi_3)+e(\chi_4)
+e(\chi_5)+e(\chi_6)\ ;$$
$$\underline z = e(\chi_0)+e(\chi_1)+e(\chi_2)+e(\chi_3)-e(\chi_4)
-e(\chi_5)-e(\chi_6)\ ;$$
$$\underline y = 6e(\chi_0)+6e(\chi_1)+6e(\chi_2)-2e(\chi_3)\ ;$$
$$\underline t = 4e(\chi_0)+4\omega^2e(\chi_1)+4\omega e(\chi_2)
-2\omega e(\chi_4)-2\omega^2 e(\chi_5)-2e(\chi_6)\ ;$$
$$\underline t^2 = 4e(\chi_0)+4\omega e(\chi_1)+4\omega^2 e(\chi_2)
-2\omega^2 e(\chi_4)-2\omega e(\chi_5)-2e(\chi_6)\ ;$$
$$\underline{tz} = 4e(\chi_0)+4\omega^2 e(\chi_1)+4\omega e(\chi_2)
+2\omega e(\chi_4)+2\omega^2 e(\chi_5)+2e(\chi_6)\ ;$$
$$\underline{t^2z} = 4e(\chi_0)+4\omega e(\chi_1)+4\omega^2 e(\chi_2)
+2\omega^2 e(\chi_4)+2\omega e(\chi_5)+2e(\chi_6)\ .$$
We show that they are all in the $\CO$-linear span of the
elements in the set 
$$\{1,z_0, z_1,s,8e(\chi_3),4e(\chi_4),4e(\chi_5)\}\ .$$
Note first that 
$$z_2 = 4\cdot 1 -z_0-z_1-s+8e(\chi_3), $$ 
$$4e(\chi_6) = 2s - 4e(\chi_4)-4e(\chi_5)$$
are in the
$\CO$-linear span of this set. One easily verifies now that
$$\underline z = 1 - s\ ,$$
$$\underline y = 6\cdot 1 - 3s - 8e(\chi_3)\ ,$$
$$\underline t = \omega z_0 + \omega^2 z_1 + z_2 
- 4\omega e(\chi_4)-4\omega^2 e(\chi_5)- 4e(\chi_6)\ ,$$
$$\underline t^2 = \omega^2 z_0 + \omega z_1 + z_2 
- 4\omega^2 e(\chi_4)-4\omega e(\chi_5)- 4e(\chi_6)\ ,$$
$$\underline{tz} = \omega z_0 + \omega^2 z_1 + z_2\ ,$$
$$\underline{t^2z} = \omega^2 z_0 + \omega z_1 + z_2\ .$$
This concludes the proof of 2.3
$\square$ \enddemo

\bigskip

The center of $\bar A = k\tenO A$ can easily be described in terms
of the generators in the quiver of $A$:

\bigskip

\proclaim{Lemma 2.4} The following set is a $k$-basis of $Z(\bar A)$.
$$\{1,\ \bar\beta\bar\gamma+\bar\gamma\bar\beta, 
\ \bar\kappa\bar\lambda+\bar\lambda\bar\kappa,\ \bar\eta\bar\delta+
\bar\delta\bar\eta,\ \bar\beta\bar\delta\bar\lambda, 
\bar\delta\bar\lambda\bar\beta,\ \bar\lambda\bar\beta\bar\delta\}\ .$$
\endproclaim

\demo{Proof} Straightforward verification, using $(\bar\beta\bar\gamma)^2 =
\bar\beta\bar\delta\bar\lambda$ and the similar relations for the
other elements in the given set.
$\square$ \enddemo

\proclaim{Proposition 2.5} For any primitive idempotent $e$ in $A$
we have $Z(A)e = eAe$. Moreover,

(i) the set $\{e_0, z_0e_0, z_1e_0, 4e(\chi_4)e_0\}$ is an $\CO$-basis
of $e_0Ae_0$.

(ii) the set $\{e_1,z_0e_1,z_2e_1,4e(\chi_4)e_1\}$ is an $\CO$-basis
of $e_1Ae_1$;

(iii) the set $\{e_2, z_1e_2, z_2e_2, 4e(\chi_5)e_2\}$ is an
$\CO$-basis of $e_2Ae_2$.

\endproclaim

\demo{Proof} Since $Z(A)\cong Z(\OAfourtilde)$ and $Z(\bar A)\cong
Z(k\tilde A_4)$, the canonical map $A \rightarrow \bar A$ maps
$Z(A)$ onto $Z(\bar A)$ and hence $Z(A)e$ onto $Z(\bar A)\bar e$.
By Nakayama's Lemma, it suffices to show that $Z(\bar A)\bar e
= \bar e\bar A\bar e$. Now $\dim_k(\bar e \bar A\bar e) = 4$
by the Cartan matrix, and so we have only to show that 
$\dim_k(Z(\bar A)\bar e) = 4$. By the symmetry of the quiver of $A$,
we may assume that $e$ corresponds to the vertex labelled $0$.
Then the set $\{\bar e, \bar\beta\bar\gamma, \bar\kappa\bar\lambda,
\bar\beta\bar\delta\bar\lambda\}$ is a $k$-basis of $Z(\bar A)\bar e$
by 2.4; in particular, $\dim_k(Z(\bar A)\bar e) = 4$ as required.
This shows that $eAe = Z(A)e$.

In order to prove (i), note that the set $$\{e_0, z_0e_0 z_1e_0,
se_0, 8e(\chi_3)e_0, 4e(\chi_4)e_0, 4e(\chi_5)e_0\}$$ generates
$e_0Ae_0$ as $\CO$-module, by the first statement and by 
the $\CO$-basis of $Z(A)$ described in 2.3.
Now we have
$$8e(\chi_3)e_0 = 2z_0e_0 - 4e(\chi_4)e_0\ ,$$
$$4e(\chi_5)e_0 = 2z_0e_0 - 2z_1e_0+4e(\chi_4)e_0\ ,$$
$$se_0 = (z_1-z_0+4e(\chi_4))e_0\ .$$
Thus the set given in (i) generates $e_0Ae_0$ as $\CO$-module,
and hence is a basis since the $\CO$-rank of $e_0Ae_0$ is $4$.
The same arguments show (ii), (iii). 
$\square$ \enddemo

\bigskip

\proclaim{Proposition 2.6} We can choose the generators $\beta$, $\gamma$,
$\delta$, $\eta$, $\lambda$, $\kappa$ in such a way that

(i) $\ A\gamma$ is the unique $\CO$-pure submodule of $Ae_0$ with character
$\chi_3+\chi_4$;

(ii) $\ A\lambda$ is the unique $\CO$-pure submodule of $Ae_0$ with character
$\chi_3+\chi_5$;

(iii) $\ A\eta$ is the unique $\CO$-pure submodule of $Ae_1$ with character
$\chi_3+\chi_6$;

(iv) $\ A\beta$ is the unique $\CO$-pure submodule of $Ae_1$ with character
$\chi_3+\chi_4$;

(v) $\ A\kappa$ is the unique $\CO$-pure submodule of $Ae_2$ with character
$\chi_3+\chi_5$;

(vi) $\ A\delta$ is the unique $\CO$-pure submodule of $Ae_2$ with character
$\chi_3+\chi_6$.
\endproclaim

\demo{Proof} We are going to prove (i); by the symmetry of the quiver of
$A$ one gets all other statements. Observe first that $\bar A\bar \gamma$
is the unique $5$-dimensional submodule of $Ae_0$ with composition factors
$2[S_0]$, $2[S_1]$, $[S_2]$. Indeed, the set $\{\bar\gamma, \bar\beta\bar\gamma,
\bar\eta\bar\gamma, \bar\gamma\bar\beta\bar\gamma, 
\bar\beta\bar\gamma\bar\beta\bar\gamma\}$ is a $k$-basis of $\bar A\bar\gamma$,
and we have $\bar\gamma, \bar\gamma\bar\beta\bar\gamma\in \bar e_0\bar A
\bar e_0$, yielding the two composition factors isomorphic to $S_0$,
we have $\bar\beta\bar\gamma, \bar\beta\bar\gamma\bar\beta\bar\gamma \in
\bar e_1\bar A\bar e_0$, yielding the two composition factors isomorphic
to $S_1$, and finally $\bar\eta\bar\gamma\in\bar e_2\bar A\bar e_0$,
yielding the remaining composition factor isomorphic to $S_2$.
One checks that there is no other submodule with exactly these composition
factors. Now there is exactly one $\CO$-pure submodule $U$ of $Ae_0$
whose reduction modulo $J(\CO)$ has composition series $2[S_0]+2[S_1]+[S_2]$,
namely the unique $\CO$-pure submodule of $Ae_0$ with character $\chi_3+
\chi_4$; this is a direct consequence of the decomposition matrix. 
One constructs $U$ as follows: write $K\tenO Ae_0 = X_0\oplus
X_3\oplus X_4\oplus X_5$, where $X_j$ is the unique submodule of
$K\tenO Ae_0$ with character $\chi_j$ for $j\in\{0,3,4,5\}$, and then
$U = Ae_0\cap (X_3\oplus X_4)$. Take now for $\gamma$ any inverse image
in $U$ of $\bar\gamma$. Then $A\gamma \subseteq U$ and $U\subseteq
A\gamma + J(\CO)U$. Thus $A\gamma = U$ by Nakayama's Lemma.
$\square$ \enddemo

\bigskip

\proclaim{Corollary 2.7} If the generators $\beta$, $\gamma$, $\delta$
$\eta$, $\lambda$, $\kappa$ are chosen such that they fulfill the
conclusions of 2.6 then, with the notation of 2.3, the following hold.

(i) $y_0\delta = y_0\eta = 0$.

(ii) $y_1\lambda = y_1\kappa = 0$.

(iii) $y_2\gamma = y_2\beta = 0$.
\endproclaim

\bigskip

\bigskip

\proclaim{Proposition 2.8} 
We can choose the generators $\beta$, $\gamma$,
$\delta$, $\eta$, $\lambda$, $\kappa$  such that the following
holds:

$$\beta\gamma = z_0e_0 = 4e(\chi_3)e_0 +  2e(\chi_4)e_1; $$
$$\gamma \beta= z_0e_1 = 4e(\chi_3)e_1 +  2e(\chi_4)e_1; $$
$$\delta \eta = z_2e_1 = 4e(\chi_3)e_1 +  2e(\chi_6)e_1; $$
$$\eta\delta = z_2e_2 = 4e(\chi_3)e_2 +  2e(\chi_6)e_2; $$
$$\lambda\kappa = z_1e_2 = 4e(\chi_3)e_2 +  2e(\chi_5)e_2; $$
$$\kappa\lambda = z_1e_0 = 4e(\chi_3)e_0 +  2e(\chi_5)e_0; $$
$$\beta\delta \lambda =  \kappa\eta\gamma= 8e(\chi_3)e_0 ; $$
$$\delta \lambda \beta= \gamma\kappa\eta = 8e(\chi_3)e_1 ; $$
$$\lambda\beta\delta = \eta\gamma \kappa = 8e(\chi_3)e_2 . $$

\endproclaim

\demo{Proof} In view of the decomposition matrix of $A$ we have
$e_0 = e(\chi_0)e_0 + e(\chi_3)e_0 + e(\chi_4)e_0 + e(\chi_5)e_0$.
Moreover, the elements $e(\chi_0)e_0$, $e(\chi_3)e_0$, $e(\chi_4)e_0$,
$e(\chi_5)e_0$ are $K$-linearly independent because they are
pairwise orthogonal idempotents in $K\tenO A$. Similar statements
hold for $e_1$, $e_2$.

We assume a choice of generators fulfilling 2.6.
We have $A\beta\gamma\subseteq A\gamma$, and the
submodule $A\gamma$ of $Ae_0$ has character $\chi_3+\chi_4$ by
2.6. Thus $\beta\gamma$ is a $K$-linear combination of
$e(\chi_3)e_1$ and $e(\chi_4)e_1$. But also $\beta\gamma$ is an
$\CO$-linear combination of the basis elements $e_1$, $z_0e_1$
$z_1e_1$, $4e(\chi_4)e_1$ given in 2.5
in which none of $\chi_1$, $\chi_5$ shows
up. Therefore $\beta\gamma$ is in fact an 
$\CO$-linear combination of the elements $z_0e_0$, $4e(\chi_4)e_0$; say
$$\beta\gamma = (\mu_0z_0e_0 +4 \nu_0 e(\chi_4))e_0 = (4\mu_0 e(\chi_3)+
2(\mu_0+2\nu_0)e(\chi_4))e_0$$ 
for some coefficients $\mu_0$, $\nu_0 \in \CO$. Hence 
$$(\beta \gamma)^2 =
(16\mu_0^2 e(\chi_3) + 4(\mu_0+2\nu_0)^2 e(\chi_4))e_0\ .$$ 
Now $(\bar\beta\bar\gamma)^2\neq 0$,
and therefore $\mu_0\in\CO^\times$. Set now
$$a_0 = 1 + \nu_0\mu_0^{-1}y_0\ .$$
Since $y_0\in J(Z(A))$ by 2.3 we have $a_0\in Z(A)^\times$.
A trivial verification, comparing coefficients, shows that we have
$$\beta\gamma = \mu_0z_0a_0e_0\ .$$
Since $\gamma = e_1\gamma = \gamma e_0$,
multiplying this with $\gamma$ on the left yields
$$\gamma\beta\gamma = \mu_0z_0a_0e_1\gamma\ .$$
Now both $\gamma\beta$ and $\mu_0z_0a_0e_1$ are contained in the
pure submodule $A\beta$ of $Ae_1$ with character $\chi_3+\chi_4$, by 
2.6 and the nature of the element $z_0$. Right multiplication by
$\gamma$ on this submodule is therefore injective (the annihilator of
$\gamma$ in $Ae_1$ is the pure submodule with character $\chi_1+\chi_6$).
Hence the previous equality implies also the equality
$$\gamma\beta = \mu_0z_0a_0e_1\ .$$
In an entirely analogous way one finds scalars $\mu_1$, $\mu_2 \in
\CO^\times$ such that, setting $a_1 = 1 + \nu_1\mu_1^{-1}y_1$ and
$a_2 = 1 + \nu_2\mu_2^{-1}y_2$, one gets the equalities
$$\delta\eta = \mu_2z_2a_2e_1\ ,\ \eta\delta = \mu_2z_2a_2e_2\ ,$$
$$\lambda\kappa = \mu_1z_1a_1e_2\ ,\ \kappa\lambda = \mu_1z_1a_1e_0\ .$$
Moreover, the equalities in 2.7 imply the following equalities:
$$a_0\delta = \delta\ , \ a_0\eta = \eta\ ,$$
$$a_1\lambda = \lambda\ , \ a_1\kappa = \kappa\ ,$$
$$a_2\gamma = \gamma\ , \ a_2\beta = \beta\ .$$
If we replace now $\beta$ by $a_0\beta$, this is not going to change
the properties stated in 2.6 and also this is not changing the relations
over $k$ of the quiver. Similarly, we can replace $\delta$ by $a_2\delta$
and $\lambda$ by $a_1\lambda$. Then the generators $\beta$, $\gamma$,
$\delta$, $\eta$, $\lambda$, $\kappa$ still fulfill 2.6, and in addition,
we have now the following equalities:
$$\beta\gamma = \mu_0z_0e_0\ ,\ \gamma\beta = \mu_0z_0e_1\ ,$$
$$\delta\eta = \mu_2z_2e_1\ ,\ \eta\delta = \mu_2z_2e_2\ ,$$
$$\lambda\kappa = \mu_1z_1e_2\ ,\ \kappa\lambda = \mu_1z_1e_0\ .$$
We have to get rid of the scalars $\mu_0$, $\mu_1$, $\mu_2$.
Since $\chi_3$ is the only character appearing in the characters of all
projective indecomposable $A$-modules we have
$$\beta\delta\lambda = 8\mu e(\chi_3)e_0$$
for some $\mu\in\CO$. Then actually $\mu\in\CO^\times$ because
$\bar\beta\bar\delta\bar\lambda \neq 0$. Moreover, $\beta\delta\lambda
\beta = 8\mu e(\chi_3)\beta$, and hence also
$$\delta\lambda\beta = 8\mu e(\chi_3)e_1\ .$$
The same argument applied again yields
$$\lambda\beta\delta = 8\mu e(\chi_3)e_2\ .$$
Applying this argument to the arrows in the quiver in the
opposite direction implies that there is $\mu'\in\CO^\times$
such that
$$\kappa\eta\gamma = 8\mu' e(\chi_3)e_0\ ,$$
$$\eta\gamma\kappa = 8\mu' e(\chi_3)e_2\ ,$$
$$\gamma\kappa\eta = 8\mu' e(\chi_3)e_1\ .$$
Now $\bar\beta\bar\delta\bar\lambda = \bar\kappa\bar\lambda
\bar\kappa\bar\lambda = \bar\kappa\bar\eta\bar\gamma$, and hence
$\mu' = \mu(1+\nu)$ for some $\nu\in J(\CO)$. Note that we can always
multiply any of the generators by any scalar in $1+J(\CO)$ without
modifying the relations over $k$. Thus, if we replace $\kappa$
by $(1+\nu)\kappa$, we may assume that $\mu' = \mu$.

Since the set $\{\kappa, \kappa\lambda\kappa\}$ is an $\CO$-basis
of $e_0Ae_2$, we can write
$$\beta\delta = a\kappa + b\kappa\lambda\kappa$$
for some unique scalars $a, b\in \CO$. Multiplying this by $\lambda$
yields
$$8\mu e(\chi_3)e_0 = \beta\delta\lambda = a\kappa\lambda + 
b(\kappa\lambda)^2 = (a\mu_1z_1 + b\mu_1^2z_1^2)e_0\ .$$
By comparing the coefficients at $e(\chi_3)e_0$ and $e(\chi_5)e_0$
of the left and right expression in this equality, we get the equations
$$8\mu = 4a\mu_1 + 16b\mu_1^2\ ,$$
$$0 = 2a\mu_1 + 4b\mu_1^2\ .$$
An easy computation shows that $b = \frac{\mu}{\mu_1^2}$.
Moreover, since $\bar\beta\bar\delta\bar\lambda = (\bar\kappa
\bar\lambda)^2$ we have $\bar a = 0$ and $\bar b = 1_k$,
hence $b = \frac{\mu}{\mu_1^2} \in 1+J(\CO)$. 
By repeating the same argument we find also that the
coefficients $\frac{\mu}{\mu_0^2}$, $\frac{\mu}{\mu_2^2}$
are in $1 + J(\CO)$. 

Next, we compute $\beta\delta\lambda\kappa\eta\gamma$ in two
different ways: on one hand we have
$$(\beta\delta\lambda)(\kappa\eta\gamma) = 64\mu^2 e(\chi_3)e_0\ ,$$
and on the other hand we have
$$\beta(\delta(\lambda\kappa)\eta)\gamma = \mu_0\mu_1\mu_2z_0z_1z_2
e(\chi_3)e_0 = 64\mu_0\mu_1\mu_2 e(\chi_3)e_0\ .$$
Together we get
$$\mu^2 = \mu_0\mu_1\mu_2\ .$$
Thus $\frac{\mu}{\mu_0^2}\frac{\mu}{\mu_1^2} = \frac{\mu_2}{
\mu_0\mu_1} \in 1 + J(\CO)$. Similarly, $\frac{\mu_1}{\mu_0\mu_2},
\frac{\mu_0}{\mu_1\mu_2} \in 1+J(\CO)$. But then also
$\frac{\mu_1\mu_2}{\mu_0}\frac{\mu_1}{\mu_0\mu_2} = \frac{
\mu_1^2}{\mu_0^2} \in 1+J(\CO)$. Since $2\in J(\CO)$ this implies
that $\frac{\mu_1}{\mu_0}\in 1+J(\CO)$. But then actually
$\mu_2 = \frac{\mu_1\mu_2}{\mu_0}\frac{\mu_0}{\mu_1} \in
1 + J(\CO)$. Similarly, $\mu_0$, $\mu_1 \in 1+J(\CO)$. 
So we can replace $\beta$ by $\mu_0^{-1}\beta$, or equivalently,
we can assume that $\mu_0 = 1$. Similarly, we can assume
that $\mu_1 = \mu_2 = 1$. Then $\mu^2 = 1$. If $\mu = -1$
we multiply all generators by $-1$; since $2 \in J(\CO)$, this does
not change the relations over $k$, but it does change the
sign of any of the above expressions $\beta\delta\lambda$
etc. involving three generators. Therefore, we can also assume
that $\mu = 1$.  

$\square$ \enddemo

We can now prove Theorem C from the introduction.

\demo {Proof of Theorem C}   We assume  a choice of generators of $A$ 
fulfilling  Proposition 2.8. 
We show that  $A$ satisfies  the  relations given in Theorem 
C. Those in the  first three  lines are obvious.
Since the set $\{\kappa, \kappa\lambda\kappa\}$ is an $\CO$-basis
of $e_0Ae_2$, we can write
$$\beta\delta = a\kappa + b\kappa\lambda\kappa$$
for some unique scalars $a, b\in \CO$. Multiplying this by $\lambda$
yields
$$8 e(\chi_3)e_0 = \beta\delta\lambda = a\kappa\lambda + 
b(\kappa\lambda)^2 =   (4a +16b)e(\chi_3) e_0 + (2a +4b)e(\chi_5) e_0  \ .$$
By comparing the coefficients at $e(\chi_3)e_0$ and $e(\chi_5)e_0$
of the left and right expression in this equality, we get the equations
$$8 = 4a + 16b\ ,$$
$$0 = 2a + 4b\ .$$
Thus the coefficients $a$, $b$ have values
$$a = -2\ ,\ b = 1\ ,$$
and from this we get  the following  relation in the statement of Theorem C:
$$\beta\delta = -2\kappa + \kappa\lambda\kappa\ .$$
In exactly the same way we get the following five relations
in the Theorem:
$$\eta\gamma = -2\lambda + \lambda\kappa\lambda\ ,$$
$$\delta\lambda = -2\gamma + \gamma\beta\gamma\ ,$$
$$\kappa\eta = -2\beta + \beta\gamma\beta\ ,$$
$$\lambda\beta = -2\eta + \eta\delta\eta\ ,$$
$$\gamma\kappa = -2\delta + \delta\eta\delta\ .$$
A similar technique is going to yield the remaining relations:
write $\gamma\beta\delta = c\delta + d\delta\eta\delta$ for
some unique $c, d \in \CO$; as before, this is possible since
$\{\delta,\delta\eta\delta\}$ is an $\CO$-basis of $e_1Ae_2$.
Multiplying by $\eta$ yields
$$\gamma\beta\delta\eta = c\delta\eta + d(\delta\eta)^2 =
cz_2e_1 + dz^2_2e_1\ .$$
The left side is equal to $(\gamma\beta)(\delta\eta) =
z_0z_2e_1$, so comparing coefficients yields now
$$16 = 4c + 16d\ ,$$
$$0 = 2c + 4d\ ,$$
and this implies $c = -4$ and $d = 2$. Thus we get indeed
$$\gamma\beta\delta = -4\delta + 2\delta\eta\delta\ $$
as claimed. The  remaining relations of this type follow in exactly the 
same way.

Now consider the  last three relations. Write 
$\beta\delta\lambda \beta =  r\beta + s\beta\gamma \beta $, for 
$ r, s \in \CO$.
Then $\beta\delta\lambda \beta \gamma =  r\beta\gamma + 
s\beta\gamma \beta\gamma $.
So
$$ 32e(\chi_3)e_0=  (4r +16s)e(\chi_3)e_0 + (2r + 4s)e(\chi_4)e_0 $$
which yields  $s=4 $ and $r=-8 $. 
The remaining two relations follow in exactly 
the same way. Thus $A$  satisfies all relations   given in Theorem C. 

Let $\tilde A$ be the 
$\CO$-algebra  described  by the generators and relations  given in 
Theorem C.   There is a surjective algebra morphism from 
$\tilde A$ to $A$. 
In order to show that $\tilde A$  and $A$ are isomorphic it suffices 
therefore to show  that  the cardinality of a minimal generating set for $A$
as an $\CO$-module  is  at most 
$24$.  Thus it suffices to check  that   the set
$$ \eqalign {\Cal {S} := &\{e_0,e_1, e_2, \beta, \gamma, \delta, \eta, \lambda, \kappa,  \cr
& \beta\gamma, \gamma\beta,  \delta\eta, \eta\delta,  \lambda\kappa, \kappa\lambda,  \cr
& \beta\gamma\beta, \gamma\beta\gamma, \delta\eta\delta, \eta\delta\eta, \lambda\kappa\lambda, \kappa\lambda\kappa, \cr
& \beta\delta\lambda, \delta\lambda \beta, \lambda\beta\delta \} \cr
}$$
 spans $\tilde A$ as $\CO$-module. This is  an easy consequence of the given 
relations; we  give  some details  for the convenience of the reader: 
 Let
$$ \Cal{G} =
\{ e_0, e_1, e_2, \beta, \gamma, \delta, \eta, \lambda, \kappa \} 
$$

From the given relations it is immediate that   for any two elements 
$x$, $y$ of $\Cal {G}$, $xy$ is  in the $\CO$-span of $\Cal {S}$. 
Thus it suffices to show that for any two elements $x, y $ of 
$\Cal {G}- \{  e_0, e_1, e_2\}  $ and any element $u$ of 
$\Cal {S} - \{ e_0,e_1, e_2, \beta, \gamma, \delta, \eta, 
\lambda, \kappa \} $,  $ xu$ and $uy$  are in the  $\CO$-span of $\Cal {S}$.   
From the  given  relations we may  also  assume that $u $ is one of 
$\beta\gamma\beta, \gamma\beta\gamma, \delta\eta\delta, \eta\delta\eta, 
\lambda\kappa\lambda, \kappa\lambda\kappa $ or  one of  
$\beta\delta\lambda, \delta\lambda \beta, \lambda\beta\delta $.

First, note that the relations $\kappa \eta  =-2\beta 
+\beta\gamma\beta  $ and 
$ \delta\lambda= -2\gamma + \gamma \beta \gamma $ give that 
$\kappa\eta\gamma = \beta\delta\lambda $. Similarly, we get 
 $\eta\gamma\kappa =\lambda \beta\delta $ and  
$\gamma \kappa\eta= \delta\lambda\beta $. 

Now suppose  $ u= \beta\gamma\beta$. Then  we may assume that $x$ is one of $ \gamma $ or 
$ \lambda $ and  that 
$y $ is one of $\gamma $ or $\delta $. The relation 
$\kappa \eta  = -2\beta +\beta\gamma \beta  $ gives 
$\gamma \kappa \eta  = -2\gamma \beta + \gamma \beta\gamma \beta  $, hence
$ \gamma \beta\gamma \beta  $  is in the $\CO$-span of $\Cal {S}$. 
The relation $\kappa\eta= -2\beta +\beta\gamma\beta $  also gives 
 $\lambda\kappa\eta= -2\lambda \beta +\lambda \beta\gamma\beta $. It follows 
from the relation  $ \lambda\kappa\eta= -4\eta + 2\eta\delta\eta $ that  
$\lambda \beta\gamma\beta $ is in the  $\CO$-span of $\Cal {S}$.  We show 
similarly that  $ \beta\gamma\beta\gamma $ and $\beta\gamma\beta \delta $ 
are in the  $\CO$-span of $\Cal {S}$.

The cases $u= \gamma\beta\gamma, \delta\eta\delta, \eta\delta\eta, \lambda\kappa\lambda, \kappa\lambda\kappa$ are handled  analogously.

Now suppose $u= \beta\delta\lambda$. Then  we may assume  that  $x$ is  
one of $\lambda $  or $\gamma $   
and $y$ is  one of $\beta $ or $\kappa$. The relation $\lambda\beta\delta\lambda = -8\lambda + 4\lambda\kappa\lambda $  shows that   
$\lambda\beta\delta\lambda$ is in   the $\CO$-span of $\Cal {S}$.  From the relation 
$\gamma \beta\delta= -4\delta  + 2\delta\eta\delta $ we get 
$\gamma \beta\delta \lambda = -4\delta\lambda  + 2\delta\eta\delta\lambda $. 
From $\gamma\kappa =-2\delta + \delta\eta\delta $, we get 
$\delta\eta\delta\lambda = \gamma \kappa \lambda +2\delta\lambda $.  Hence 
$\delta\eta\delta\lambda $ is in   the $\CO$-span of $\Cal {S}$, and  so is 
$\gamma \beta\delta \lambda$. We argue similarly  to show that  
$\beta\delta\lambda\beta $ and $ \beta\delta\lambda\kappa $ are in   the 
$\CO$-span of $\Cal {S}$.

The cases $u= \delta\lambda \beta$ and $u=\lambda \beta\delta $ are handled   in the same fashion.
$\square$\enddemo

\bigskip

\noindent {\bf Remark 2.9.} 
An interesting consequence of 2.5 is the
structure of $eAe$ for any primitive idempotent $e$ in $A$.
We have an $\CO$-algebra isomorphism
$$eAe \cong \CO[X,Y]/<X^2-Y^2-2(X-Y)\ ,\ XY-2X^2+4X>\ ;$$
indeed, we may assume that $e = e_0$, and then the assignment
$X\mapsto z_0e_0$, $Y\mapsto z_1e_0$ induces the required
isomorphism. In particular, we have an isomorphism of
$k$-algebras
$$\bar e\bar A \bar e \cong k[X,Y]/<X^2-Y^2\ ,\ XY>\ .$$
This is, by Erdmann [\Erd, III.1, III.3], up to isomorphism the unique $4$-dimensional
symmetric $k$-algebra which is not isomorphic to the group algebra
of the Klein four group.
One might be tempted to ask whether any symmetric $\CO$-algebra is
the endomorphism algebra of some projective module of some block
algebra.

\bigskip

%\noindent {\bf Remark 2.10.} 
%Similar computations in the remaining case where $kG\bar b$ is Morita equivalent
%to the principal block algebra of $k\tilde A_5$ lead to showing that
%$\OG b$ is Morita equivalent to an algebra $A$ which can be described in
%terms of its quiver and relations; however, certain scalars in those relations
%remain undetermined and it is thus still an open problem whether $A$ must
%be Morita equivalent to the principal block algebra of $\CO \tilde A_5$. 

\bigskip\bigskip\bigskip
\enddocument